\documentclass[12pt]{amsart}
\usepackage{amscd,amssymb}
\usepackage[graph,frame,poly,arc]{xy}  
\usepackage[plainpages,backref,urlcolor=blue]{hyperref}
\usepackage{mathtools}

\theoremstyle{plain}
\newtheorem{thm}[subsection]{Theorem}
\newtheorem{lem}[subsection]{Lemma}

\newtheorem{cor}[subsection]{Corollary}

\theoremstyle{definition}
\newtheorem{rk}[subsection]{Remark}
\newtheorem{definition}[subsection]{Definition}

\newtheorem{conj}[subsection]{Conjecture}
\newtheorem{ass}[subsection]{Assumption}

\numberwithin{equation}{section}
\newcommand{\I}{{\mathcal I}}

\newcommand{\A}{{\mathcal A}}
\newcommand{\B}{{\mathcal B}}

\newcommand{\CC}{{\mathcal C}}

\newcommand{\HH}{{\mathcal H}}

\newcommand{\al}{{\alpha}}

\newcommand{\SSS}{{\mathcal S}}

\newcommand{\C}{\mathbb{C}}
\newcommand{\PP}{\mathbb{P}}

\newcommand{\dd}{{\rm d}}

\DeclareMathOperator{\defect}{def}

\DeclareMathOperator{\reg}{reg}

\begin{document}

\title [Pascal's Theorem and  Ziegler's line arrangements]
{From Pascal's Theorem to the geometry of Ziegler's line arrangements}

\author[Alexandru Dimca]{Alexandru Dimca}
\address{Universit\'e C\^ ote d'Azur, CNRS, LJAD, France and Simion Stoilow Institute of Mathematics,
P.O. Box 1-764, RO-014700 Bucharest, Romania}
\email{dimca@unice.fr}

\author[Gabriel Sticlaru]{Gabriel Sticlaru}
\address{Faculty of Mathematics and Informatics,
Ovidius University
Bd. Mamaia 124, 900527 Constanta,
Romania}
\email{gabriel.sticlaru@gmail.com }



\begin{abstract} 
G\" unter Ziegler has shown in 1989 that some homological invariants associated with the free resolutions of Jacobian ideals of line arrangements are not determined by combinatorics. His classical example involves hexagons inscribed in conics. Independently, Sergey Yuzvinsky has arrived in 1993 at the same type of line arrangements in order to show that formality is not determined by the combinatorics.
In this note we look into the geometry of such line arrangements, and find out an unexpected relation to the classical Pascal's Theorem. Our results give information on the minimal degree of a Jacobian syzygy and on the formality of such hexagonal line arrangements in general, without an explicit choice for the six vertices of the hexagon.
\end{abstract}
 
\maketitle


\section{Introduction} 
Let $S=\C[x,y,z]$ be the polynomial ring in three variables $x,y,z$ with complex coefficients, and let $C:f=0$ be a reduced curve of degree $d\geq 3$ in the complex projective plane $\PP^2$. 
We denote by $J_f$ the Jacobian ideal of $f$, i.e. the homogeneous ideal in $S$ spanned by the partial derivatives $f_x,f_y,f_z$ of $f$, and  by $M(f)=S/J_f$ the corresponding graded quotient ring, called the Jacobian (or Milnor) algebra of $f$.
Consider the graded $S$-module of Jacobian syzygies of $f$ or, equivalently, the module of derivations killing $f$, namely
\begin{equation}
\label{eqD0}
D_0(f)=\{\rho=(a,b,c) \in S^3 \ : \ af_x+bf_y+cf_z=0\}.
\end{equation}

We say that $C:f=0$ is an {\it $m$-syzygy curve} if the minimal number of generators of the module  $D_0(f)$ is $m$. Then the module $D_0(f)$ is generated by $m$ homogeneous syzygies, say $\rho_1,\rho_2,...,\rho_m$, of degrees $d_j=\deg \rho_j$ ordered such that $$1 \leq d_1\leq d_2 \leq ...\leq d_m.$$ 
We call these degrees $(d_1, \ldots, d_m)$ the {\it exponents} of the curve $C$ and $\rho_1,...,\rho_m$ a {\it minimal set of generators } for the module  $D_0(f)$. 
The smallest degree $d_1$ is sometimes denoted by ${\rm mdr}(f)$ and it is called the minimal degree of a Jacobian relation for $f$. 

The curve $C$ is {\it free} when $m=2$, since then  $D_0(f)$ is a free module of rank 2, see for instance \cite{KS,ST,To}.
We say that Terao's Conjecture  holds  for a free line arrangement $\A$ if any other line arrangement $\B$  having an  isomorphic intersection lattice  $L(\B)=L(\A)$,  is also free, see \cite{DHA,OT,Yo14}. 
Terao's Conjecture is known to hold in many cases, for instance when the number of lines in $\A$ is at most 14, see \cite{BK}. One of the  main conjectures in the area is then the following. For a more general version, refer to \cite[Conjecture 4.138]{OT}.
\begin{conj}
\label{conj0}
Terao's Conjecture  holds  for any free line arrangement $\A$.
\end{conj}

Since the total Tjurina number $\tau(\A)$ is determined by the intersection lattice $L(\A)$,  a possible approach to proving Terao's Conjecture may be  to check that $\A:f=0$ and $\B:g=0$ satisfy ${\rm mdr}(f)={\rm mdr}(g)$ and then apply the characterization of free curves 
as being the curves with maximal total Tjurina number for a fixed degree
$d$ and $r={\rm mdr}(f) <d/2$, see \cite{Dmax, duPCTC}. Hence Conjecture \ref{conj0} would be a consequence of the following stronger claim.

\begin{conj}
\label{conj1}
For any  line arrangement $\A:f=0$, the integer $r={\rm mdr}(f)$
is determined by the intersection lattice when $r<d/2$.
\end{conj}
In view of \cite[Proposition 3.2]{AD}, this Conjecture is equivalent to a positive answer to
the question asked in \cite[Question 7.12]{Cook+}.

The additional condition $r<d/2$ in Conjecture \ref{conj1} comes from examples due to G. Ziegler in \cite[Example 8.7]{Z}, see also Remark \ref{rkZ10} below, which show that 
the invariant ${\rm mdr}(f)$ is not combinatorially determined in general. Indeed, Ziegler produced a pair of line arrangements $\A_Z:f_Z=0$ and $ \A_Z': f_{Z}'=0$ containing both 9 lines,
6 triple points and 18 double points, with the same intersection lattice
$L(\A_Z)=L(\A_Z')$ and  such that 
$$5={\rm mdr}(f_Z)< {\rm mdr}(f_Z')=6.$$
A similar pair of arrangements $(\A,\A')$ was discovered by S. Yuzvinski, see \cite[Example 2.2]{Y}, where it is shown that $\A'$ is 2-formal, while $\A$ is not. Recall that for a line arrangement $\A$ consisting of $d$ lines $L_j:\ell_j=0$ for $j=1, \ldots, d$, we define the {\it relation space } $F(\A)$ to be the kernel of the 
evaluation map
\begin{equation}
\label{2form}
\C^d \to S_1, \  e_j  \mapsto \ell_j,
\end{equation}
where $e_1,...,e_d$ is the canonical basis of $\C^d$.
We say that $\A$ is 2-{\it formal} if $F(\A)$ is spanned by relations involving only 3 lines in $\A$, see \cite[Definition 1.2]{SchT}. The 2-formality is simply called {\it formality} in \cite[Definition 2.1]{DST}.
Yuzvinsky's example is also discussed from this point of view in \cite[Example 1.4]{SchT}.
The pair of arrangements $(\A_Z,\A_Z')$ and $(\A,\A')$ were obtained as follows: start with 6 points
$A_1, \ldots , A_6$ in the plane $\PP^2$, consider the 6 sides of the associated hexagon $\HH$, namely the lines $L_j: \ell_j=0$ determined by the points
$A_j$ and $A_{j+1}$, for $j=1, \ldots, 6$, where we set $A_7=A_1$.
Add the 3 diagonals $L'_k: \ell_k'=0 $ of the hexagon $\HH$ determined by the points
$A_k$ and $A_{k+3}$ for $k=1,2,3$. Assume that the following assumption holds.
\begin{ass}
\label{Ass0}
The resulting line arrangement 
$$\A(A_1, \ldots,A_6)=L_1 \cup L_2 \cup \ldots \cup L_6 \cup L'_1 \cup L'_2 \cup L'_3$$
consists of 9 distinct lines and has only double points except the 6 triple points situated at the points $A_1, \ldots , A_6$. 
\end{ass}
Then the Ziegler arrangement
$\A_Z$ corresponds to one choice of the points $A_1, \ldots , A_6$, where they
are situated on the (degenerate) conic
\begin{equation}
\label{QZ1}
Q_Z:(x-y-z)(2x-y+z)=0
\end{equation}
and the Ziegler arrangement
$\A_Z'$ corresponds to another choice of the six points $A_1, \ldots , A_6$, where they
are not situated on a conic. The relation with conics was in fact noticed later, by H. Schenck, see \cite[Example 13]{Sch}, where the author concentrates for some unknown reasons on {\it smooth} conics.
For more details concerning this arrangement $\A_Z$ see Section 3 below. 
Other such pairs of arrangements have been constructed, by varying the choice of the points $A_1, \ldots , A_6$, see for instance the arrangements $\A_D$ and $ \A'_D$ in  \cite[Remark 8.5]{DHA}. Here for the arrangement $\A_D$ the points $A_1, \ldots , A_6$ are chosen on the smooth conic
\begin{equation}
\label{Q1}
Q:x^2+y^2-z^2=0,
\end{equation}
see Section 4 below for more details on this arrangement $\A_D$.

There is an {\it interesting question to find new pairs of line arrangements},
say $\B:g=0$ and $\B':g'=0$ such that their intersection lattices verify $L(\B)=L(\B') \ne L(\A_Z)$ and
${\rm mdr}(g)\ne {\rm mdr}(g').$ One can construct such pairs $\B$ and $\B'$ by adding suitable lines to the Ziegler type arrangement pairs
$\A_Z, \A_Z'$ or $\A_D, \A_D'$, see \cite[Remark 3.3]{OP}. In particular, in Remark \ref{rkZ10} below we construct such a {\it new pair of line arrangements} of degree $d=10$  by adding a line to the arrangements $\A_Z$ and $\A'_Z$ such that their intersection lattices verify $L(\B)=L(\B')$ and
\begin{equation}
\label{eqZP10}
5={\rm mdr}(g)= \frac{d}{2} <{\rm mdr}(g')=6.
\end{equation}
This explains the strict inequality $r<d/2$ in Conjecture \ref{conj1}.
A different and much deeper approach for the construction of Ziegler type pairs  $\B$ and $\B'$, at least in the case $r>d/2$,  can be found in \cite{DST}.

In this note we study first the geometry of the line arrangements $\A_Z$ and $\A _D$ and try to understand why they have this special behavior. To state our results, first we recall some more notation. For any reduced plane curve $C:f=0$,
let $J=J_f$ be the Jacobian ideal of $f$ and $I=I_f$ be its saturation with respect to the maximal ideal $(x,y,z)$. Then the singular subscheme $\Sigma_C$ of the reduced curve $C:f=0$ is the 0-dimensional scheme defined by the ideal $I$ and we define the following sequence of defects
\begin{equation} 
\label{eqDEF}
\defect _k\Sigma_C=\tau(C)-\dim \frac{S_k}{I_k}.
\end{equation}
Here $\tau(C)$ is the total Tjurina number of $C$, that is the sum of all Tjurina numbers of the singularities of $C$.
With this notation, one has the following result, a consequence of \cite[Theorem 1]{Bull13}.

\begin{thm}
\label{thm1}
Let $C:f=0$ be a reduced degree $d$ curve in $\PP^2$ with $r= {\rm mdr }(f) <d-1$. Then
$$\defect _k\Sigma_C=\tau(C)-\dim \frac{S_k}{I_k}=0$$
for $k > 2d-5-r$  and
$$\defect _{k}\Sigma_C=\tau(C)-\dim \frac{S_{k}}{I_{k}}>0$$
for $k \leq 2d-5-r$.
In particular 
$$\defect _k\Sigma_{\A_Z}=\defect _k\Sigma_{\A_ D}=0$$
for $k>8$ and 
$$\defect _8\Sigma_{\A_Z}=\defect _8\Sigma_{\A _D}>0.$$
\end{thm}
Since $\tau(\A_Z)=\tau(\A_D)=42$, $\dim S_8=45$ and $\dim J_8=3$
for both arrangements $C=\A_Z$ and $C=\A_D$, it follows that there is at least one octic form $h_C \in I_8 \setminus J_8$. 
Then $h_C \notin J_8$ means that the form $h_C$ is not
in the net spanned by $f_x,f_y,f_z$. And $h_C \in I_8$ means that for each singular point $p$ of $C$, the germ of the regular function $h_{C,p}$ associated to $h_C$ at $p$
is in the local ideal $\I_p$ defining the germ $(\Sigma_C,p)$ of the singular subscheme $\Sigma _C$ at $p$.
\begin{definition}
\label{defSAT}
We  say in the sequel simply that $h_C$ belongs to $\I_p$ when $h_{C,p} \in \I_p$.
\end{definition}
Our first main result in Theorem \ref{PascalThm} says that this octic $h_C$ has a nice, unexpected geometry.
Let $B_k=L_k \cap L_{k+3}$ be the intersection of the {\it opposite sides} in the hexagon $\HH$ and recall that {\it Pascal's Theorem} tells us that the points $B_1,B_2$ and $B_3$ are situated on a line, call it $L: \ell=0$, when
the vertices of $\HH$ are situated on a conic, see \cite{GH}, p. 673. This can be interpreted as saying that the linear form $\ell$ belongs to the ideals
$\I_p$ for $p \in \{B_1,B_2,B_3\}$. This line $L$ plays a key role in our story.

Let $A_j'$ be the intersection of the lines $L_j$ and $L_{j+2}$ for $j=1, \ldots,6$,
where we set $A_7=A_1$ and $A_8=A_2$. 

\begin{thm}
\label{PascalThm}
Assume that the vertices $A_1, \ldots, A_6$ of the hexagon $\HH$ are situated on a conic, which is either smooth or a union of two distinct lines. Consider the associated line arrangement
$\A=\A(A_1, \ldots, A_6)$
 and assume that $\A$ satisfies Assumption \ref{Ass0}. Then there is a  nonzero
quartic form $q_{\A} \in S_4$ such that $q_{\A}$ vanishes on the set of 12 points
 $$ \{A_1,A_2,A_3,A_4,A_5,A_6,A_1',A_2',A_3',A_4',A_5',A_6' \}$$
 and the corresponding octic form
$$h_{\A}=\ell \cdot \ell_1'\cdot \ell_2'\cdot \ell_3'\cdot q_{\A}$$
belongs to $I_8$ but not to $J_8$.
\end{thm}
In view of Theorem \ref{thm1}, our Theorem \ref{PascalThm} yields the following.
\begin{cor}
\label{PascalCor}
Assume that the vertices $A_1, \ldots, A_6$ of the hexagon $\HH$ are situated on a conic, which is either smooth or a union of two distinct lines. Consider the associated line arrangement
$\A=\A(A_1, \ldots, A_6)$
 and assume that $\A$ satisfies Assumption \ref{Ass0}.
Then, if we set $\A:f=0$, one has 
${\rm mdr}(f) =5$ and $\A$ is not 2-formal.
\end{cor}

The case when the vertices $A_1, \ldots, A_6$ of the hexagon $\HH$ are not situated on a conic is settled in our second main  result.
\begin{thm}
\label{thm2}
Assume that the vertices $A_1, \ldots, A_6$ of the hexagon $\HH$ are not situated on a conic. Consider the associated line arrangement
$\A=\A(A_1, \ldots, A_6)$
 and assume that $\A$ satisfies Assumption \ref{Ass0}.
Then, if we set $\A:f=0$, one has $I_8=J_8$, 
${\rm mdr}(f) =6$ and $\A$ is 2-formal.
\end{thm}
The proofs of Theorem \ref{PascalThm} and Theorem \ref{thm2} are given in Section 2 and  use
the following approach. To test the fact that a form $h$ belongs to the saturation $I$ of the Jacobian ideal $J$, we use the localization criterion described before Definition \ref{defSAT}. Each double point $p$ of $\A$ imposes a linear condition on $h$, and each triple point $q$ of $\A$ imposes 4 linear conditions. Indeed, the 1-jet of $h$ at $q$ has to vanish (3 conditions), and there is one extra condition on the 2-jet of $h$ at $q$, see Lemmas \ref{lemD4} and \ref{lemD4bis} below. 
In the case of Theorem \ref{PascalThm} this approach produces a system $\SSS$ of 
18 linear equations involving 15 unknowns, and we have to show that
its rank is $<15$. By using the fact that the nodes of a nodal curve impose independent conditions on some linear systems, we can replace this system by a new system $\SSS'$ of 6 linear equations involving 3 unknowns, and we have to show that
its rank is $<3$, which can be done either by using a computer algebra software as SINGULAR \cite{Sing}, or by hand.

In the case of Theorem \ref{thm2} this approach produces a system $\SSS$ of 
42 linear equations involving 45 unknowns, and we have to show that
its rank is 42. 
Using some  results due to  H. Schenck in \cite{Sch2} and to R. Burity, S. O. Toh\u aneanu and Yu Xie in \cite{BTX}, we replace this system by a new system $\SSS'$ of 6 linear equations involving 9 unknowns, and we have to show that
its rank is $6$. This can again be done either by using a computer algebra software as SINGULAR or by hand.

 The claims about 2-formality follow using a recent key result of  M. DiPasquale, J. Sidman and W. Traves, see
\cite[Corollary 7.8]{DST}.

\medskip

We would like to thank S. Toh\u aneanu for very useful discussions, in particular concerning the 2-formality of line arrangements.

\section{Proof of the main results } 

\subsection{Proof of Theorem \ref{PascalThm}}

Consider the line arrangement $\B$ determined by the 6 sides of the hexagon $\HH$. Then $\B$ has exactly 15 double points, 12 of them located at the 12 points
$A_j$ and $A'_j$ for $j=1, \ldots,6$ and 3 of them at the points $B_1$, $B_2$ and $B_3$. Using \cite[Corollary 1.6]{Edin}, it follows that vanishing at these 15 points imposes 15 independent conditions on the polynomials in $S_4$. Since $\dim S_4=15$, it follows that the vector space $Z_4$ of quartic forms $q$ vanishing at  the 12 points
$A_j$ and $A'_j$ for $j=1, \ldots,6$ has dimension $15-12=3$.
The following 3 quartic forms 
\begin{equation} 
\label{eqB1}
q_1=\ell_2\ell_3\ell_5\ell_6,  \  q_2=\ell_1\ell_3\ell_4\ell_6 \text{ and } q_3=\ell_1\ell_2\ell_4\ell_5
\end{equation}
are clearly in $Z_4$. Moreover, they are linearly independent, as any two of them has in common 2 linear factors which do not divide the third one.
It follows that any quartic form $q$ in $Z_4$ can be written as
$$q=c_1q_1+c_2q_2+c_3q_3$$
for some coefficients $c_j \in \C$. 
If we set 
$$h=\ell \ell'_1\ell'_2 \ell'_3 q$$
 for such a quartic $q$, then clearly $h$ vanishes at all the singular points of $\A$.
But the points
$A_j$ are triple points, and each such point adds a new linear condition
on the coefficients $c_j$, if we impose the condition in Definition \ref{defSAT} to hold at $p=A_j$.
To understand this condition, we need the following.
\begin{lem}
\label{lemD4}
 Let $R=\C\{u,v\}$ be the convergent power series local ring with $\C$ coefficients and variables $u$ and $v$. Let $g \in R$ be a power series such that its 3-jet is a product of 3 distinct nonzero linear factors in $u,v$, namely
  $$j^3g=\al_1 \al_2 \al_3.$$
   Let $f \in R$ be a power series such that its 2-jet is given by $j^2f=\al_1\al$, for some linear form $\al$ in $u,v$. Then
 $f$ belongs to the local Jacobian ideal $ J_g=(g_u,g_v) \subset R$ of $g$  if and only if $\al$ is proportional to the linear form which is the coefficient of $\dd u \wedge \dd v$ in the product
 of differentials
 $$\dd \al_1 \wedge \dd (\al_2 \al_3).$$
\end{lem}
\proof Note that the condition on $j^3g$  says exactly that the singularity $g=0$ has type $D_4$,
and hence it is isomorphic to a triple point of a line arrangement.
The claim is clear if we recall that for any $D_4$ singularity $g$ one has
${\bf m}_2^3 \subset J_g$, where ${\bf m}_2=(u,v)$ is the maximal ideal of $R$. This implies that  $f \in J_g$ if and only if   $j^2f \in J_{j^3g}$.
The fact that this last condition is equivalent to the claimed proportionality can be checked by a direct computation.
\endproof
To see what happens at the triple points of $\A$, let us consider the point
$p=A_1$. There are 3 lines in $\A$ passing through $p$, namely
two sides $L_1$, $L_6$ of the hexagon $\HH$ and the diagonal $L_1': \ell_1'=0$.
Note that $\ell_1'$ is a factor of $h_{\A}$ and hence when we localize at the point $p$ the condition that $h_{\A}$ belongs to $\I_p$, we are exactly in the setting of Lemma \ref{lemD4}, where the linear form $\al$ is just a local equation for the tangent line to the quartic $\CC: q_{\A}=0$ at the point $p$. Hence we have to fix the direction of all the tangent lines at the quartic $\CC$ at the points $A_j$ for $j=1, \ldots, 6$.
This gives 6 linear equations with 3 unknown $c_1,c_2,c_3$. To show that this system has a non trivial solution, we use a simple direct computation with the software SINGULAR.

\begin{figure}[ht]
\includegraphics[scale=0.35]{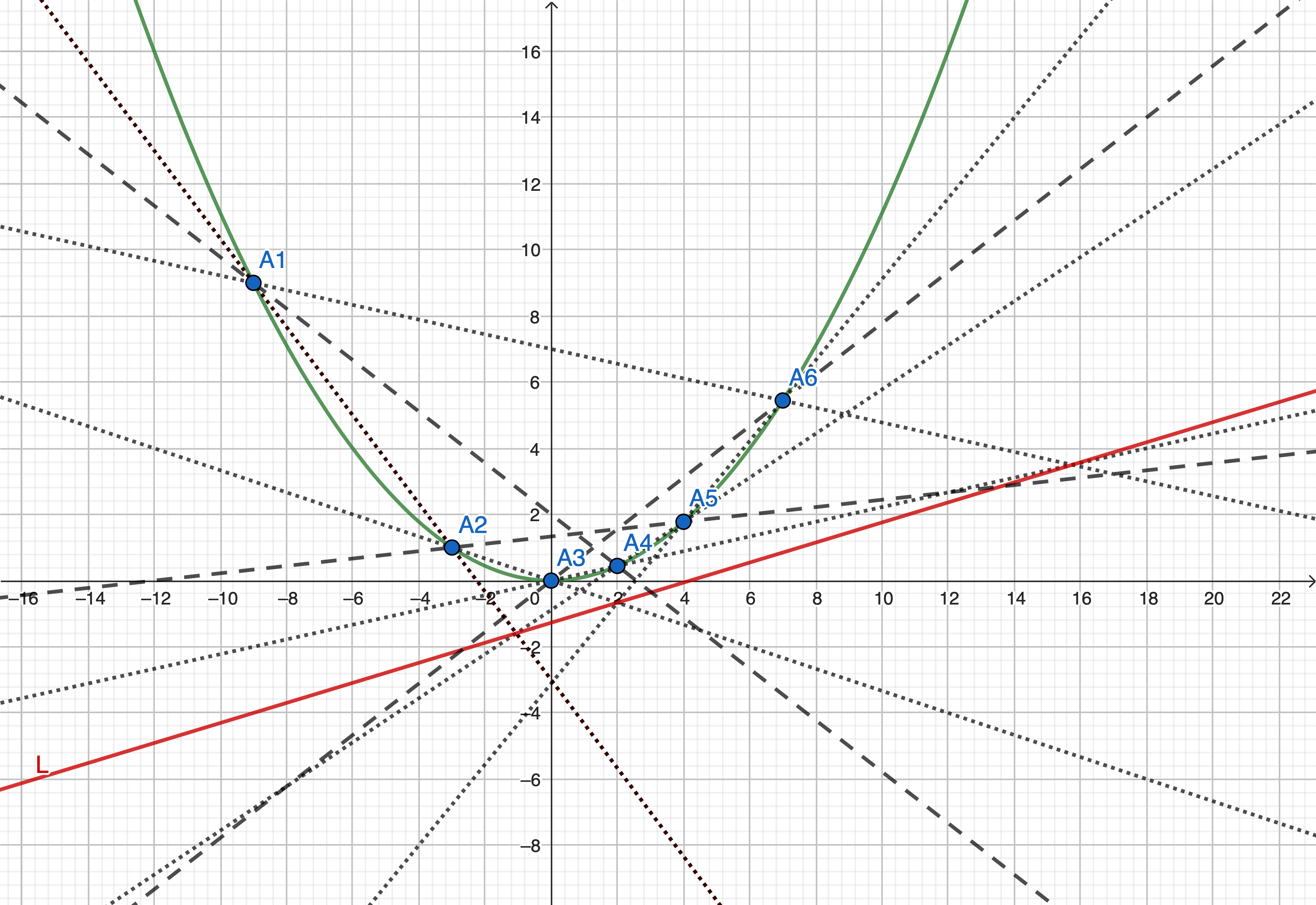}
\caption{Hexagon in a smooth conic}
\label{Parabole1}
\end{figure}

We give some more details in the case when the conic containing the points $A_j$ for $j=1, \ldots, 6$ is smooth. Then we can choose the coordinates on $\PP^2$ such that the conic
has the equation $x^2=yz$.
We assume that  $A_j=(t_j:t_j^2:1)$, with $t_j \in \C$, for $j=1, \ldots, 6$.
In the local coordinates $(u,v)$ at $A_1$ defined by $x=t_1+u$, $y=t_1^2+v$ and $z=1$, using Lemma \ref{lemD4}, we get the following formula for the tangent line
$T_1$ at the quartic $\CC$
$$T_1: (t_1t_2+t_1t_6+2t_2t_6-2t_1t_4-t_2t_4-t_4t_6)u-(t_2+t_6-2t_4)v=0$$
in order that the germ of $h$ belongs to $\I_p$ for $p=A_1$.
Writing that this equation is proportional to the usual equation for this tangent line, namely
$$T_1':h_x(t_1,t_1^1,1)u+h_y(t_1,t_1^2,1)u=0,$$
we get a first linear equation in the coefficients $c_j$.
A simple computation shows that
$$T_1': (c_1\ell_3(t_1)(t_1+t_6)+c_3\ell_4(t_1)(t_1+t_2))u-(c_1\ell_3(t_1)+c_3\ell_4(t_1)v=0.$$
Then the first equation of our linear system is, after dividing by $(t_1-t_4)(t_2-t_6)\ne 0$, the following 
\begin{equation} 
\label{eqS1}
E_1:  \  (t_1-t_3)(t_4-t_6)c_1+(t_1-t_5)(t_2-t_4)c_3=0.
\end{equation} 
The  corresponding formulas  for the other tangent directions $T_j$ at $\CC$ at $A_j$ can be obtained from the above formulas by cyclically permuting the points $t_j$, for $j=1, \ldots, 6$ and in the same time, the unknown $c_1,c_2,c_3$.
In this way we get the next 5 equations
\begin{equation} 
\label{eqS2}
E_2:  \  (t_2-t_6)(t_3-t_5)c_1+(t_2-t_4)(t_5-t_1)c_2=0,
\end{equation} 
\begin{equation} 
\label{eqS3}
E_3:  \  (t_3-t_1)(t_4-t_6)c_2+ (t_3-t_5)(t_6-t_2)c_3=0,
\end{equation} 
\begin{equation} 
\label{eqS4}
E_4: \  (t_4-t_6)(t_1-t_3)c_1+ (t_4-t_2)(t_5-t_1)c_3=0,
\end{equation} 
\begin{equation} 
\label{eqS5}
E_5:  \  (t_5-t_3)(t_6-t_2)c_1+ (t_5-t_1)(t_2-t_4)c_2=0,
\end{equation} 
\begin{equation} 
\label{eqS6}
E_6:  \  (t_6-t_4)(t_1-t_3)c_2+ (t_6-t_2)(t_3-t_5)c_3=0.
\end{equation} 
The resulting $6 \times 3$ matrix is shown to have all the 3-minors equal to 0, either by a hand computation or using SINGULAR. Hence the nonzero quartic $q=q_{\A}$ such that $h_{\A} \in I_8$ exists. Indeed,  all the nodes in $\A$ distinct from the points $A_j'$ are situated on at least one of the 4 lines $L,L_1',L_2',L_3'$.

\medskip

When the conic containing the points $A_j$ is a union of 2 distinct lines, one can choose these 2 lines to be the union of two coordinate axes $xy=0$ and the points $A_j$ as the points with the following coordinates
$$(0:t_1:1), (t_2:0,1),(0:t_3:1), (t_4:0:1),(0:t_5:1) \text{ and } (t_6:0:1)$$
in this order. Then the quartics in \eqref{eqB1} give, exactly as above, a basis for the quartics vanishing at the 12 points $A_j$ and $A'_j$ for $j=1, \ldots,6$.
The corresponding formulas for the double points
for the tangent directions $T_j$ and $T'_j$ are easy to determined, and the resulting linear system of 6 equations in $c_1,c_2,c_3$ has again rank 2, and hence admits a non trivial solution.

The fact that $h_{\A}$ is not in $J_8$ can be seen as follows. Note that the equation $f=0$ for $\A$ can be written as $f=\ell_1' w$ for some polynomial $w$. Then a linear relation
$$h_{\A}=af_x+bf_y+cf_z$$
evaluated at a general point $q$ of $L_1'$ with $w(q) \ne 0$, yields 
$$a\ell_{1,x}'+b\ell_{1,y}'+c\ell_{1,z}'=0.$$
In other words, we get $(a:b:c)\in L_1'$. Since we can repeat this argument with $L_2'$ and $L_3'$, and since these 3 lines have no common point, this gives a contradiction.

\begin{rk}
\label{rkA}

(i) Any nontrivial solution $(c_1,c_2,c_3)$ of the system above has all the components non zero. Indeed, if we assume for instance that
$c_3=0$, then equation \eqref{eqS1} implies that $c_1=0$, and the second equation \eqref{eqS2} would give $c_2=0$.
The fact that all the components $c_j$ are non zero implies that the quartic $\CC: q_{\A}=0$ is smooth at all the points $A_j$,
for $j=1, \ldots,6$. Indeed, at each point $A_j$, one of the quartics $q_k$ is singular, and the other two are smooth, with distinct tangent spaces.

(ii) Note that the  quartic curve $\CC$ passes through all triple points and through the 6 double points $A'_j$, not situated on the other components of $h_Z=0$. In view of Lemma \ref{lemD4}, the tangent to $\CC$ at each of the triple points $p$ is uniquely determined (unless $p$ is a singular point for $\CC$). It follows that $\CC$ is {\it the unique quartic curve with these properties} if $\CC$ is irreducible. Indeed, if $\CC'$ would be another such curve, the intersection multiplicity
$(\CC,\CC')_p$ would be at least 1 at each node $p$ and at least 2 at each triple point $p$. It follows that
$$(\CC,\CC') \geq 6 + 12= 18 > 16=\deg \CC \deg \CC',$$
a contradiction, unless $\CC$ is reducible.
\end{rk}

\subsection{Proof of Corollary  \ref{PascalCor}}
Theorem \ref{thm1} and Theorem \ref{PascalThm} implies that 
${\rm mdr}(f) \leq 5$. We show now that ${\rm mdr}(f) >4$, that is that
$D_0(f)_4=0$. To do this, notice that the arrangement $\A$ is obtained from the nodal arrangement $\B:g=0$ considered in the proof of Theorem \ref{PascalThm} by adding the 3 diagonals.
Since $\B$ is nodal, it follows that ${\rm mdr}(g) =4$, see \cite[Theorem 4.1]{Edin}. Let $\B':g'=0$ be the line arrangement obtained from $\B$ by adding the diagonals $L_1'$ and $L_2'$. Then one has
${\rm mdr}(g') \geq 4$, see for instance \cite{DIS}. Apply now \cite[Theorem 6.2]{DIS}, which is a consequence of the main result in \cite{STY}, for $k=4$, with the first curve $C_1$ being $\B'$ and the second curve $C_2$ being the last diagonal $L_3'$.
Note that $L_3'$ intersects $\B'$ at 6 points, namely the triple points 
$A_3$ and $A_6$ and 4 nodes. The exact sequence in 
\cite[Theorem 6.2]{DIS} implies then that $D_0(f)_4=0$.

Now we prove the claim about the 2-formality. The key result here is
\cite[Corollary 7.8]{DST}, saying that an arrangement $\A:f=0$ consisting of $d$ lines which is irreducible is
not 2-formal if and only if
\begin{equation} 
\label{eqF1}
\reg (D_0(f))=d-2,
\end{equation} 
where $\reg(N)$ denotes the Castelnuovo-Mumford regularity of a graded $S$-module $N$ and $\A$ irreducible means $\A$ has only points of multiplicity $\leq d-2$. As noticed in the proof of \cite[Corollary 3.5]{DIM}, one has the equality
\begin{equation} 
\label{eqF2}
\reg (D_0(f))=\reg(M(f))-d+3,
\end{equation} 
where $M(f)=S/J_f$ denotes the Jacobian or Milnor algebra of $f$.
Moreover, \cite[Theorem 3.3]{DIM} says that one has
\begin{equation} 
\label{eqF3}
\reg(M(f))=st(f)-1
\end{equation} 
for any line arrangement $\A:f=0$ which is not free. Here $st(f)$ is
 the {\it stability threshold} of $f$ defined by
$${st}(f)=\min \{q~~:~~\dim M(f)_k=\tau(\A) \text{ for all } k \geq q\}.$$
In addition, for any  arrangement $\A$ with $d \geq 4$ lines one has
\begin{equation} 
\label{eqF31}
st(f)\leq 2d-4,
\end{equation} 
see  \cite[Corollary 3.5]{DIM}.
It follows that, for non free line arrangements, the condition \eqref{eqF1} can be restated as
\begin{equation} 
\label{eqF4}
st(f)=2d-4,
\end{equation} 
in other words $\A$ is not 2-formal if and only if $st(f)$ is maximal.
Now \cite[Theorem 1]{Bull13} implies that
\begin{equation} 
\label{eqF5}
\dim M(f)_{2d-5}=\dim M(f_s)_{d-1} + \defect _{d-1} \Sigma _{\A},
\end{equation} 
where  $f_s$  a homogeneous polynomial in $S$ of the same degree $d\geq 4$ as $f$ and such that $C_s:f_s=0$ is a smooth curve in $\PP^2$.
This implies that 
$$\dim M(f_s)_{d-1}=\dim S_{d-1}-3= \binom{d+1}{2}-3.$$
If we apply \eqref{eqF5} to our line arrangement $\A$,  we have $d=9$ and
we know by Theorem \ref{PascalThm} that 
$ \defect_{8}\Sigma_{\A}>0$. Hence
we get
$$ \dim M(f)_{13}=45-3+ \defect_{8}\Sigma_{\A}>42=\tau(\A),$$
which implies $st(f) \geq 14$. In view of \eqref{eqF31} this implies $st(f) =14$ which
proves our claim by \eqref{eqF4}.

\subsection{Proof of Theorem \ref{thm2}}

Assume that the vertices $A_1, \ldots, A_6$ of the hexagon $\HH$ are not on a conic. Consider the line arrangement
$$\A=\A(A_1, \ldots, A_6)=L_1 \cup \ldots \cup L_6 \cup L_1' \cup L_2' \cup L_3'$$
satisfying Assumption \ref{Ass0}.
Then, if we set $\A:f=0$, one has
${\rm mdr}(f) \leq 6$, see \cite{Mich} or \cite[Example 2.4]{BT}. Moreover, in view of Theorem \ref{thm1}, the equality ${\rm mdr}(f) =6$ is equivalent to
$\dim I_8=\dim J_8=3$. 
Note that
$$f= \ell_1 \ell_2\ell_3 \ell_4\ell_5 \ell_6\ell_1' \ell_2'\ell_3'$$
and set 
$$g_j=f/\ell_j \text{ for } j=1,2, \ldots, 6 \text{ and } h_k=f/\ell_k'\text{ for } k=1,2,3.$$

Let $K$ be the ideal in $S$ generated by the 9 polynomials
$g_1, \ldots, g_6,h_1,h_2,h_3$. It is clear that the Jacobian ideal $J=J_f$ is contained in $K$, and hence its saturation $I=J^{sat}$ is contained in $K^{sat}$.
Then \cite[Lemma 3.2]{Sch2} or \cite[Theorem 2.2]{BTX} says that the ideal $K$ and the quotient $S/K$ have linear graded free resolutions. Next \cite[Remark 1.1]{BTX} says that in these conditions one has
$$K_m=K_m^{sat}$$
for any $m \geq 8$, since $K$ is generated in degree 8. It follows that
\begin{equation} 
\label{eqT1}
I_8 \subset K_8^{sat}=K_8.
\end{equation} 
Hence any polynomial $\phi \in I_8$ can be written as a linear combination
\begin{equation} 
\label{eqT2}
\phi =t_1g_1+t_2g_2+ \ldots +t_6g_6+s_1h_1+s_2h_2+s_3h_3.
\end{equation} 
Any such linear combination vanishes at all the double points in $\A$ and has a singularity, namely a zero 1-jet, at any triple point $A_j$ of $\A$.
We have to show that the 6 conditions $\phi \in \I_{A_j}$ for $j=1,...,6$ as in Definition \ref{defSAT}, 
impose 6 linearly independent conditions on the 9 coefficients $t_j,s_k$.
To do this, we use the following.

\begin{lem}
\label{lemD4bis}
 Let $R=\C\{u,v\}$ be the convergent power series local ring with $\C$ coefficients and variables $u$ and $v$. Let $g \in R$ be a power series such that its 3-jet is a product of 3 distinct nonzero linear factors in $u,v$, namely
  $$j^3g=\al_1 \al_2 \al_3.$$
   Let $f \in R$ be a power series such that its 2-jet is given by $$j^2f=a\al_2\al_3+b\al_3\al_2+c\al_1\al_2,$$
    for some constants $a,b,c \in \C$. Then
 $f$ belongs to the local Jacobian ideal $ J_g=(g_u,g_v) \subset R$ of $g$  if and only if $a,b,c$ satisfy the following equation
 $$a\dd \al_2\wedge \dd \al_3+b\dd \al_3\wedge \dd \al_1+c\dd \al_1\wedge \dd\al_2=0.$$
\end{lem}
\proof 
Note that any quadratic form in $u,v$ can be written as a sum
$$a\al_2\al_3+b\al_3\al_2+c\al_1\al_2.$$ 
Moreover, the equation 
 $$a\dd \al_2\wedge \dd \al_3+b\dd \al_3\wedge \dd \al_1+c\dd \al_1\wedge \dd\al_2=0$$
 defines a 2-dimensional subspace in the 3-dimensional space of quadratic form in $u,v$. The proof of this claim is exactly as the proof of
 Lemma \ref{lemD4} above.
 \endproof
 We can assume that the points $A_1$, $A_3$ and $A_5$ are not collinear. Indeed, if they are collinear, we replace them by $A_2$, $A_4$ and $A_6$. These new 3 points cannnot be collinear, due to our assumption that the vertices of the hexagon $\HH$ are not on a conic.
 Since $A_1$, $A_3$ and $A_5$ are not collinear, we can choose the coordinates on $\PP^2$ such that
$$A_1=(1:0:0), \  A_3=(0:1:0)  \text{  and  } A_5=(0:0:1).$$
Then we set
$$A_j=(x_j:y_j:z_j)  \text{  for } j =2,4,6.$$
These coordinates $(x_j:y_j:z_j)$ satisfy some obvious conditions.  For instance
\begin{equation} 
\label{eqT3}
z_6 \ne 0,
\end{equation} 
since otherwise the points $A_1$, $A_3$ and $A_6$ are collinear, and hence the side $L_6$ would coincide with the diagonal $L_3'$. Moreover
\begin{equation} 
\label{eqT4}
x_2y_6-x_6y_2 \ne 0,
\end{equation} 
since otherwise the points $A_2$, $A_5$ and $A_6$ are collinear,
and hence the side $L_5$ would coincide with the diagonal $L_2'$.
Consider the matrix
\begin{center}
$$N=\left(
  \begin{array}{ccc}
     x_2y_2 & y_2z_2& x_2z_2 \\
    x_4y_4 & y_4z_4& x_4z_4 \\
    x_6y_6 & y_6z_6& x_6z_6 \\
    \end{array}
\right)$$
\end{center}
and set $\Delta= \det (N)$. Then
\begin{equation} 
\label{eqT5}
\Delta \ne 0
\end{equation} 
is the condition that the 6 points $A_j$ are not on a conic. Using the above coordinates, it is easy to get the following formulas for the equations of our lines:
$$\ell_1=z_2y-y_2z, \ \ell_2=z_2x-x_2z, \  \ell_3= z_4x-x_4z, \ \ell_4=y_4x-x_4y, \ell_5=y_6x-x_6y,$$
$$\ell_6=z_6y-y_6z, \ \ell_1'=z_4y-y_4z, \  \ell_2'=y_2x-x_2y \text{ and } 
\ell_3'=z_6x-x_6z.$$
Now we use Lemma \ref{lemD4bis} to find the 6 equations satisfied by the 9 coefficients $t_j,s_k$ in \eqref{eqT2}. Consider first the point $A_1$. At this point, only the polynomials $g_1$, $g_6$ and $h_1$ have a non zero second order jet, hence the condition here involves only $t_1$, $t_6$ and $s_1$. Since we can write the 3-jet of $f$ at $A_1$ up to a non-zero multiplicative constant as
$$f=\ell_1\ell_6\ell_1',$$
the condition on $t_1$, $t_6$ and $s_1$ comes from the equation
$$t_1\dd \ell_6 \wedge \dd \ell_1' +t_6\dd \ell_1' \wedge \dd \ell_1+s_1\dd \ell_1 \wedge \dd \ell_6=0.$$
This yields the equation
\begin{equation} 
\label{eqT6}
E_1: (y_6z_4-y_4z_6)t_1+(y_4z_2-y_2z_4)t_6+(y_2z_6-y_6z_2)s_1=0.
\end{equation} 
Similar computations at the points $A_3$ and $A_5$ yield the following
two equations.
\begin{equation} 
\label{eqT7}
E_3: (x_6z_4-x_4z_6)t_2+(x_2z_6-x_6z_2)t_3+(x_4z_2-x_2z_4)s_3=0.
\end{equation}
\begin{equation} 
\label{eqT8}
E_5: (x_2y_6-x_6y_2)t_4+(x_4y_2-x_2y_4)t_5+(x_6y_4-x_4y_6)s_2=0.
\end{equation}
When we compute the equation corresponding to the point $A_2$ using the same approach, we get 3 forms, namely
$\dd x \wedge \dd y$, $\dd x \wedge \dd z$ and $\dd y \wedge \dd z$. Indeed, one can check that Lemma \ref{lemD4bis} can be reformulated such that localization at $A_2$ is not necessary, one can work in the homogeneous coordinates $(x,y,z)$ as well. The key point is that the linear forms $\al_1, \al_2, \al_3$ in the coordinates $(x,y,z)$ must have a common point $p$, and the computation is valid at this point $p$.
The coefficients of these 3 forms have all the form
$(x_2t_1-y_2t_2+z_2s_2)A$, where the factor $A$ takes the values
$z_2,y_2,x_2$, hence we have $A \ne 0$ as in \eqref{eqT3} above.
It follows that the equation corresponding to the point $A_2$ is the following
\begin{equation} 
\label{eqT9}
E_2:  x_2t_1-y_2t_2+z_2s_2=0.
\end{equation}
Similar computations at the points $A_4$ and $A_6$ yield the following
two equations.
\begin{equation} 
\label{eqT10}
E_4: -y_4t_3+z_4t_4+x_4s_1=0.
\end{equation}
\begin{equation} 
\label{eqT11}
E_6: z_6t_5+x_6t_6-y_6s_3=0.
\end{equation}
Consider now the matrix $M$ with 6 rows and 9 columns obtained as follows: the first 3 rows come from the simpler equations $E_2$, $E_4$ and $E_6$, and the last 3 rows come from the equations $E_1$, $E_3$ and
$E_5$. It remains to show that this matrix has rank 6, that is that there is a 6-minor $m_6\ne 0$. Using SINGULAR,  we get that the minor corresponding to the columns  $1,2,3,4,5,7$ has the equation
$$m_6=-x_4x_6^2y_2^2y_6z_2z_4z_6+x_2x_6^2y_2y_4y_6z_2z_4z_6+x_2x_4x_6y_2y_6^2z_2z_4z_6-x_2^2x_6y_4y_6^2z_2z_4z_6+$$
$$+x_4x_6^2y_2^2y_4z_2z_6^2-2x_2x_4x_6y_2y_4y_6z_2z_6^2+x_2^2x_4y_4y_6^2z_2z_6^2-x_2x_6^2y_2^2y_4z_4z_6^2+$$
$$+x_2x_4x_6y_2^2y_6z_4z_6^2+x_2^2x_6y_2y_4y_6z_4z_6^2-x_2^2x_4y_2y_6^2z_4z_6^2.$$
A direct computation with SINGULAR shows that
\begin{equation} 
\label{eqT12}
m_6=z_6(x_2y_6-x_6y_2)\Delta,
\end{equation}
and hence $m_6 \ne 0$ in view of \eqref{eqT3}, \eqref{eqT4} and \eqref{eqT5}.
This completes the proof of the first two claims in Theorem \ref{thm2}.

The claim concerning the 2-formality follows using the same approach as in the proof of Corollary  \ref{PascalCor} above. The key difference is that now
$ \defect_{8}\Sigma_{\A}=0$ since $\dim I_8=\dim J_8=3$. Hence
we get
$$ \dim M(f)_{13}=42+ \defect_{8}\Sigma_{\A}=42=\tau(\A),$$
which says that $st(f)\leq 13<2d-4=14$. In view of \eqref{eqF4}, this  proves our claim.

\section{The octic $h_{\A_Z}$ and a Ziegler pair of degree 10} 

In this section we study as an example the octic $h_{\A}$ in the case of Ziegler's line arrangement  
\begin{equation}
\label{Z1}
\A _Z: f_Z=xyz(x+y-z)(x-y+z)(2x-2y+z)(2x-y-2z)
\end{equation}
$$ (2x+y+z)(2x-y-z)=0, $$
introduced  \cite{Z}.
Using SINGULAR, we see that the corresponding saturated ideal $I$ has 5 generators, four of degree 8 and one of degree 9. The four generators of degree 8 have proportional remainders modulo the Jacobian ideal $J$, and hence this common remainder (up to a non-zero constant factor) can be used as
the additional element to be added to $J_8$ to get $I_8$. Computation with SINGULAR shows that this element can be taken to be 
\begin{equation}
\label{hZ}
 h_Z(x,y,z)=xyz(y-3z)q_Z(x,y,z),
\end{equation}
where 
\begin{equation}
\label{qZ}
q_Z(x,y,z)=12x^4-12x^3y+3x^2y^2-6xy^3+3y^4-4x^3z+12x^2yz-2xy^2z+3y^3z-
\end{equation}
$$-13x^2z^2+6xyz^2-6y^2z^2+2xz^3-3yz^3+3z^4.$$
Using $GeoGebra$ we can draw the following picture.

\begin{figure}[ht]
\includegraphics[scale=1]{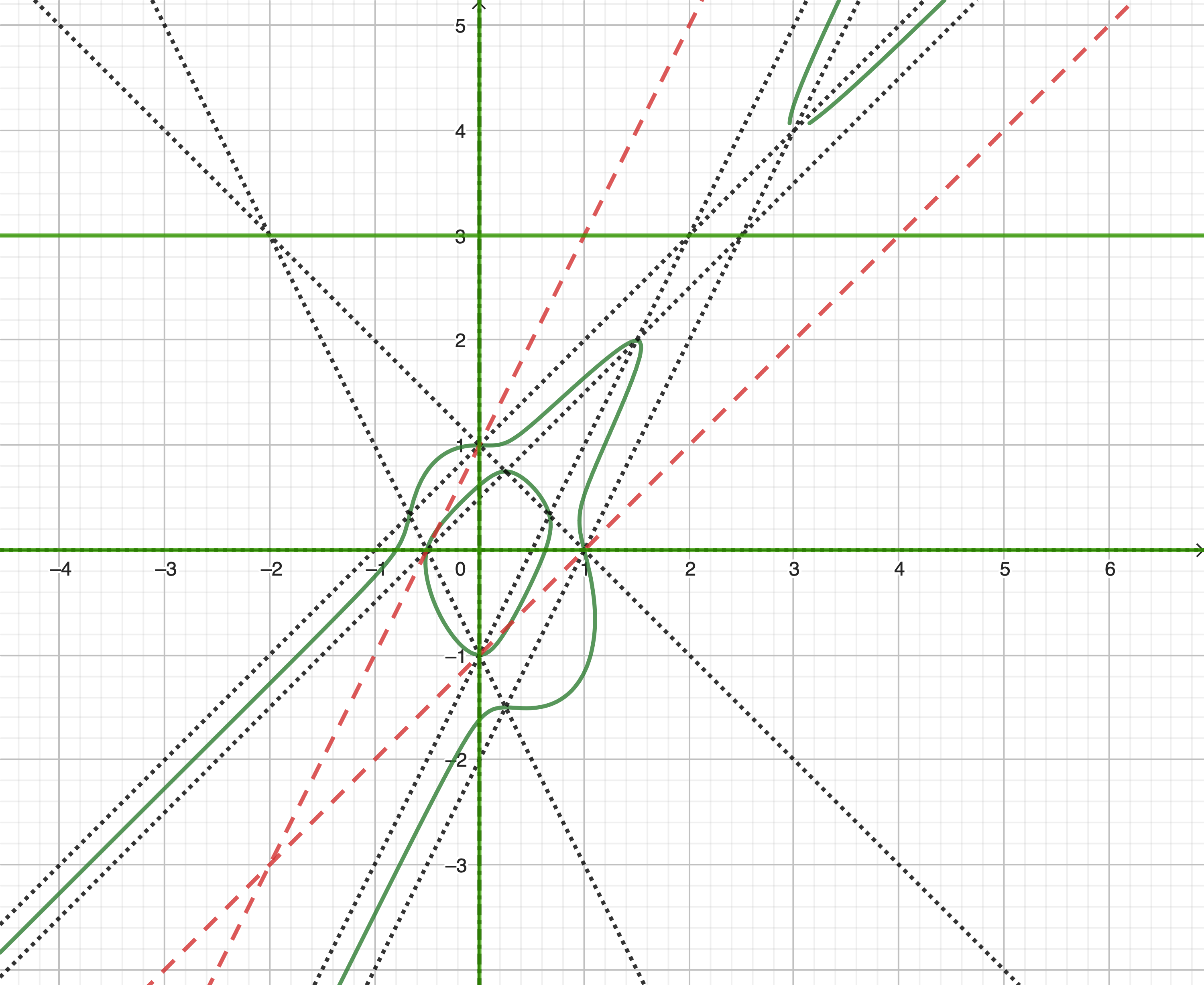}
\caption{Ziegler's arrangement $\A_ Z$}
\label{Ziegler1}

\end{figure}
The lines in the arrangements $\A_ Z$ are dotted (black) lines, with the line at infinity $z=0$ not drawn. The degenerate conic $Q_Z$ from \eqref{QZ1}
is represented by the two dashed (red) lines.
The 5 factors of $h_Z$ correspond to 

\begin{enumerate}
\item[(i)] The 3 lines $x=0$, $y=0$ and $z=0$ in the arrangements $\A _Z$, which are the diagonals of the corresponding hexagon $\HH$.
\item[(ii)] The line $L:y-3z=0$ determined by the 3 double points, corresponding to the intersection of the opposite sides in the hexagon $\HH$ formed by the 6 lines in $\A _Z$, not listed at point $(i)$. These 3 double points are collinear in view of Pascal's Theorem.
\item[(iii)] A quartic curve $\CC_Z$, whose set of real points, drawn in a continuous (green) line has 2 connected topological components in the corresponding projective plane.
\end{enumerate}

\begin{rk}
\label{rkZ10}
If we add to the line arrangement $\A_Z$ the line $L_Z:x-y-z=0$, which is one of the two lines forming the degenerate conic $Q_Z$ in \eqref{QZ1},
computation with SINGULAR shows that the obtained line arrangement
$\B_Z=\A_Z \cup L_Z: g=0$ satisfies ${\rm mdr}(g)=5=d/2$, where
$d= \deg g=10$. The line arrangement $\A_Z':f'=0$ can be taken to be as the line arrangement obtained from $\A_Z$ by moving the triple point
$p=(0:1:1)$ to the new position $p'=(0:2:1)$, hence outside the conic $Q_Z$. Define $\B_Z'=\A_Z '\cup L_Z: g'=0$. Since $\B_Z'$ has degree
$d=10$ and has points of multiplicity $4$ coming from the 3 triple points of $\A_Z$ situated on the line $L_Z$, it follows from \cite[Theorem 1.2]{Mich}
that one has ${\rm mdr}(g') \leq 6$. This implies that   ${\rm mdr}(g')= 6$ as we have claimed in Introduction, see \eqref{eqZP10}.
Indeed, it is well known, see for instance \cite{DIS}, that by addition of a line to a curve, the invariant ${\rm mdr}$ cannot decrease, hence
$${\rm mdr}(g') \geq {\rm mdr}(f')=6.$$
\end{rk} 

\section{The octic $h_{\A_D}$} 

Consider the following two arrangements 
\begin{equation}
\label{myZ1}
\A_D: f_D=xy(x-y-z)(x-y+z)(2x+y-2z)(x+3y-3z)(3x+2y+3z)
\end{equation}
$$(x+5y+5z)(7x-4y-z)=0,$$
and respectively by
\begin{equation}
\label{myZ2}
\A': f'=xy(4x-5y-5z)(x-y+z)(16x+13y-20z)(x+3y-3z)(3x+2y+3z)
\end{equation}
$$(x+5y+5z)(7x-4y-z)=0$$
see  \cite[Remark 8.5]{DHA}, but beware a misprint in the equation for $\A'$ given there.
This pair of arrangements  satisfy ${\rm mdr}(f)=5$ and ${\rm mdr}(f')=6$, though $\A$ and $\A '$ have the same combinatorics. 

In this section we study the octic $h_{\A}$ for the line arrangement  $\A_D:f_D=0$. The corresponding saturated ideal $I$ has again 5 generators, four of degree 8 and one of degree 9. The four generators of degree 8 have again proportional remainders modulo the Jacobian ideal $J$, and hence this common remainder (up to a non-zero constant factor) can be used as
the additional element to be added to $J_8$ to get $I_8$. Computation with SINGULAR shows that this element can be taken to be 
\begin{equation}
\label{h}
 h(x,y,z)=xy(7x-4y-z)(x-y-19z)q(x,y,z),
\end{equation}
where 
\begin{equation}
\label{q}
q(x,y,z)=51x^4+8x^3y+29x^2y^2+17xy^3+75y^4+12x^3z-13x^2yz+
\end{equation}
$$ +19xy^2z-15y^3z-186x^2z^2+ +26xyz^2-210y^2z^2-12xz^3+15yz^3+135z^4   .$$
Using $GeoGebra$ we can draw the following picture.

\begin{figure}[ht]
\includegraphics[scale=0.6]{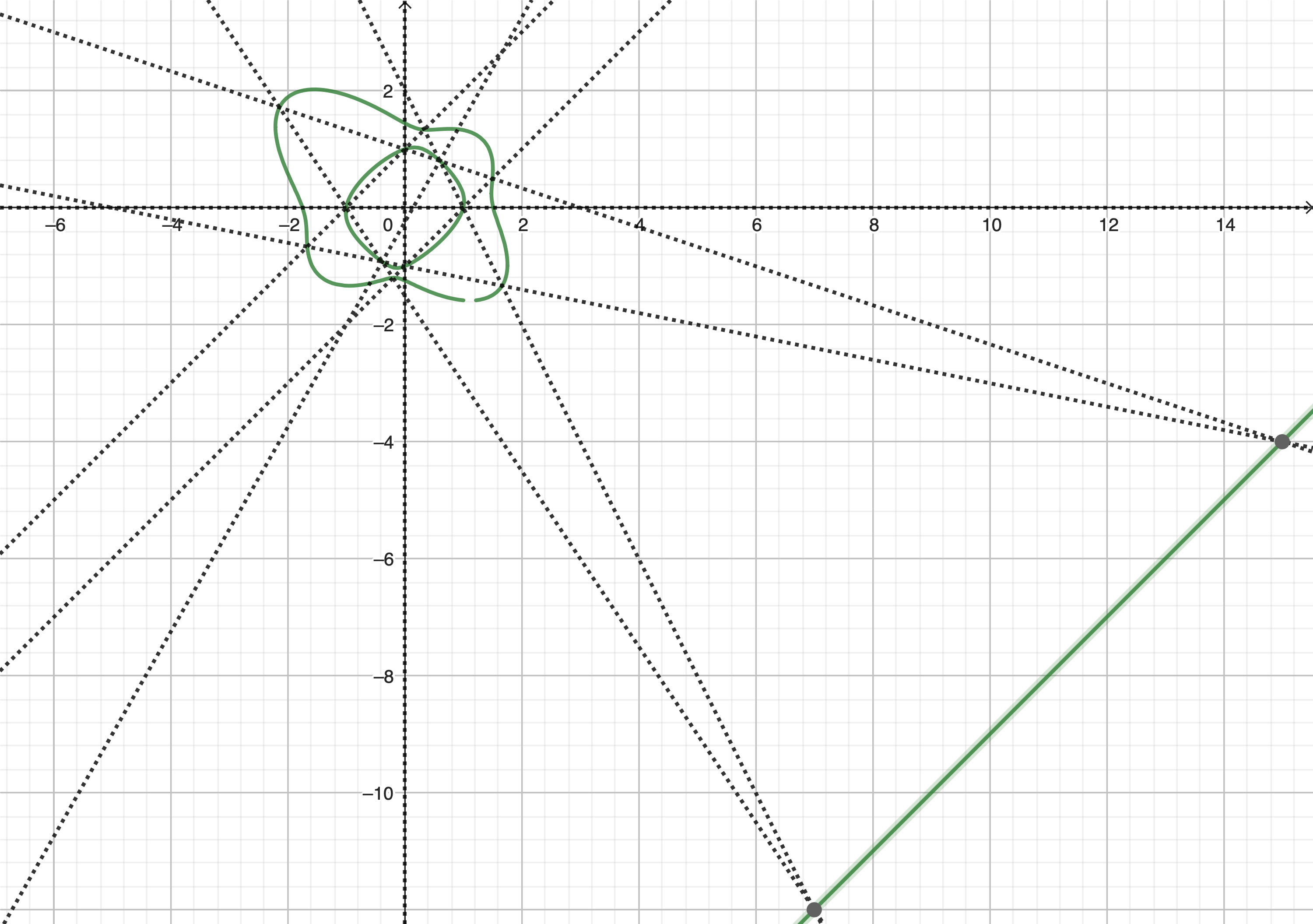}
\caption{Line arrangement $\A_D$}
\label{MyZiegler1}
\end{figure}

The 9 lines in the arrangements $\A$ are dotted (black) lines. The smooth conic $Q$ from \eqref{Q1} is the unit circle $x^2+y^2-1=0$,
but it is not represented in the picture, to keep this picture as clear as possible.
The 5 factors of $h$ correspond to 

\begin{enumerate}
\item[(i)] The 3 lines $x=0$, $y=0$ and $7x-4y-z=0$ in the arrangements $\A $, which are exactly  the diagonals of the hexagon $\HH$ formed by the other 6 lines in $\A$.
\item[(ii)] The line $L: x-y-19z=0$, draw in a continuous (green) line, is determined by the 3 double points, corresponding to the intersection of the opposite sides in the hexagon $\HH$. These 3 double points are collinear again in view of Pascal's Theorem.
\item[(iii)] A quartic curve $\CC$, whose set of real points, drawn in a continuous (green) line has 2 connected topological components.
\end{enumerate}

\begin{rk}
\label{rkZD}
The polynomial $f_Z$ has smaller coefficients than the polynomial $f_D$,
and this makes all the computations easier with this polynomial. On the other hand, the line arrangement $\A_Z$ is more difficult to picture in the affine plane given by $z=1$, since one of the line is the line at infinity
$z=0$, which contains two of the triple points. The reader may compare Figure \ref{Ziegler1} to Figure \ref{MyZiegler1} to understand this point.
\end{rk}

\end{document}